\documentclass[letterpaper, 10 pt, conference]{ieeeconf} 
\IEEEoverridecommandlockouts                              
\usepackage{cite}
\usepackage{amsmath,amssymb,amsfonts}
\usepackage{algorithmic}
\usepackage{graphicx}
\usepackage{textcomp}
\usepackage{xcolor}
\usepackage{amsmath,amssymb,amsfonts,mathrsfs}
\newtheorem{assumption}{Assumption}

\newtheorem{remark}{Remark}
\newtheorem{theorem}{Theorem}
\newtheorem{lemma}{Lemma}

\newtheorem{definition}{Definition}

\begin{document}

\title{Adaptive linear quadratic optimal control for discrete-time stochastic systems\\
\author{Nian Liu and Cheng Zhao
\thanks{This work was supported by . }
\thanks{N. Liu and C. Zhao is with the State Key Laboratory of Mathematical Sciences,  Academy of Mathematics and Systems Science, Chinese Academy of Sciences, Beijing 100190, P. R. China;.
(e-mail: liunian@amss.ac.cn, zhaocheng@amss.ac.cn.)}}
}

\maketitle

\begin{abstract}
In this paper, we consider the adaptive linear quadratic optimal control problem for a class of discrete-time stochastic systems, where the system matrix pair $(A,B)$ is unknown and the state weighting matrix $Q$ in the cost function is not necessarily positive semidefinite. Our proposed adaptive optimal control only assumes that the corresponding algebraic Riccati equation has a stabilizing solution, which weakens the standard controllable and observable assumptions used in the literature. To address this problem, we modify a weighted least squares algorithm by incorporating a random regularization method, which ensures the uniform stabilizability of the family of estimated models. Additionally, a diminishing excitation is introduced in the design of the proposed adaptive control to guarantee the strong consistency of the desired components of the estimates. Finally, by leveraging this family of estimates and employing a suitable switching method, we demonstrate that a certainty equivalence control with such a diminishing excitation  is an optimal adaptive control. 
\end{abstract}


\section{Introduction}
Linear quadratic (LQ) optimal control is one of the most important issues and has been thoroughly studied in the field of system and control. For example, it is well-known that the optimal control can be easily obtained for the infinite-horizon LQ control problem, provided that the corresponding algebraic Riccati equation (ARE) admits a stabilizing solution. However, due to the widespread presence of uncertainties in system structure, information, and the environment, it is natural to consider adaptive optimal control where the system matrices are assumed to be unknown. 

The adaptive linear quadratic control problem, where the system matrices are unknown, presents a key theoretical challenge to ensure both the global stability of the closed-loop adaptive systems and the controllability of the estimated models \cite{aa1,aa2,aa3,a8,a1,a2}. This issue was first addressed for discrete-time open-loop stable systems (autoregressive moving average with exogenous input models)  in \cite{a8}. Later, in \cite{a1}, the author removed the assumption of open-loop stability, by proving the self-convergence property of the weighted least squares (WLS) algorithm. The adaptive continuous-time linear quadratic Gaussian control problem, under the assumption that $(A, B)$ is controllable and $(A,Q^{1/2})$ is observable, was later effectively addressed in \cite{a2}, utilizing the self-convergence property of the continuous-time WLS algorithm and a random regularization theory established in \cite{a1}. To the best of our knowledge, the assumptions used in the adaptive LQ optimal control problem are strictly stronger than those in the LQ optimal control problem with known system parameters. A natural question to ask is whether it is always possible to design an adaptive optimal control when there exists an optimal control in the corresponding LQ control problem. In order to deal with this issue, in \cite{g2}, we construct an optimal adaptive control under the weaker assumption that $(A, B)$ is stabilizable and $(A, Q^{1/2})$ is detectable. However, such condition is still stronger (and not necessary) than the condition that the corresponding ARE admits a stabilizing solution. 

In this paper, we consider the adaptive linear quadratic control problem for a class of discrete-time stochastic system and give a concrete approach for the design of the optimal adaptive control, under the assumption that corresponding algebraic Riccati equation has a stabilizing solution. Our proposed optimal adaptive control is designed through the following steps. The initial step involves employing the WLS algorithm to generate a family of convergent estimates, without imposing constraints on excitation or closed-loop system stability. In the second step, to obtain uniformly stabilizable estimated models \cite{g2}, the estimates of the system matrices are refined through a random regularization procedure outlined in \cite{a1}. Finally, by integrating a diminishing excitation (white noise) and certainty equivalence principle into the switching control, we demonstrate that the desired components of the estimates of system matrices are strongly consistent, and at the same time, our proposed adaptive control is optimal for the quadratic cost function. 

The main contributions of this paper are summarized as follows: 
\begin{enumerate}
    \item We extend the recent continuous-time adaptive LQ control results under partial parameter identifiability in \cite{g2} to the discrete-time setting.
    \item We prove that adaptive optimal control exists if and only if the corresponding known-parameter LQ problem admits a stabilizing Riccati solution, completely closing the gap between the existence theories of adaptive and known-parameter optimal control.
    \item We show that full parameter consistency is not required for adaptive optimality, revealing a fundamental separation between parameter identification and control performance.
\end{enumerate}

The remainder of the paper is organized as follows. In Section \ref{section2}, we introduce the discrete-time adaptive LQ control problem and the WLS estimator associated with the random regularization method. In Section \ref{section3}, the design procedure of the optimal adaptive control is presented and the main results on the stability and
optimality of the closed-loop systems are established. Section \ref{section4} gives the proof of the main theorems. Finally, some concluding remarks are given in Section \ref{section6}.

\section{Problem Formulation}\label{section2}
\subsection{LQ optimal control}

Consider the following discrete-time LQ control systems:
\begin{equation}\label{s1}
	x(t+1)=Ax(t)+Bu(t)+Dw(t),
\end{equation}
where $x(t)\in \mathbb{R}^n$ is the state, $u(t)\in\mathbb{R}^{m}$ is the control, $A\in\mathbb{R}^{n\times n}, B\in\mathbb{R}^{n\times m},D\in \mathbb{R}^{n\times p}$ are \emph{unknown} matrices, and $\{w(t)\in \mathbb{R}^p,\mathscr{F}_t\}$ is a marginal difference sequence defined on the basic probability space $(\Omega,\mathscr{F},P)$ with
$$\lim\limits_{T\to \infty}\frac{1}{T}\sum_{t=1}^{T} w(t)w^{\mathsf{T}}(t)= I_p\quad a.s.$$
and $$ \sup_{t\ge0}E[|w(t+1)|^2|\mathscr{F}_t]<\infty\quad a.s.$$

The objective is to design an admissible control to minimize the following ergodic cost function:
\begin{equation}\label{cost}
	J(u)=\limsup\limits_{T\to \infty}\frac{1}{T}\sum_{t=1}^{T}(x^{\mathsf{T}}(t)Qx(t)+u^{\mathsf{T}}(t)Ru(t)),
\end{equation}
where the symmetric matrix  $Q$ and positive definite matrix $R$ are known weighting matrices.  

First, we introduce the definition of admissible controls.
\begin{definition}
	A control $\{u(t),t\ge0\}$ for the system (\ref{s1}) is said to be admissible if the following property holds for any initial state $x(0)\in\mathbb{R}^n$:
	\begin{equation}\notag
 \limsup\limits_{T\to \infty}\frac{\sum_{t=0}^{T}|u(t)|^2 }{\sum_{t=0}^{T}|x(t)|^2 }<\infty\quad a.s.
	\end{equation}
\end{definition}
\vskip 0.2cm

It is well-known that if the following ARE:
\begin{equation}\label{are1}
    X=Q+A^{\mathsf{T}}XA-A^{\mathsf{T}}XB(R+B^{\mathsf{T}}XB)^{-1}B^{\mathsf{T}}XA
\end{equation}
or equivalently
$$X=A^{\mathsf{T}}X(I_n+BR^{-1}B^{\mathsf{T}}X)^{-1}A+Q$$
admits a stabilizing solution $X$, i.e., $$A-B(R+B^{\mathsf{T}}XB)^{-1}B^{\mathsf{T}}XA$$ is stable, then the feedback control
\begin{equation}\label{u}
    u(t)=-(R+B^{\mathsf{T}}XB)^{-1}B^{\mathsf{T}}XAx(t)
\end{equation}
can stabilize the system and minimize the cost function. Moreover, the minimal cost is
\begin{equation}
    J(u)=tr(D^{\mathsf{T}}XD) \quad a.s.
\end{equation}

Since the system matrices $A,B$ are unknown, the optimal feedback control (\ref{u}) is not available. Therefore, it is natural to consider the adaptive control. It is worth mentioning that, the papers \cite{a1} and \cite{a2} design the adaptive optimal control  for continuous-time and discrete-time systems respectively, under the conditions that $(A,B)$ is controllable and $(A,Q^{1/2})$ is observable. The paper \cite{g2} considers the continuous-time LQG problem under the stabilizability and detectability assumptions. In this paper, we will show that the following weaker condition is sufficient for the system (\ref{s1})-(\ref{cost}) to design an adaptive optimal control. 

\begin{assumption}\label{ass1} 
The ARE (\ref{are1}) has a stabilizing solution.
\end{assumption}
\begin{remark}
    In \cite{p13}, it is shown that the ARE (\ref{are1}) has a stabilizing solution if and only if the pair $(A,B)$ is stabilizable and the matrix function $G(s)$ is antianalytic perfactorizable, where
	$G(s)=R+B^T(\frac{1}{s}I-A^T)^{-1}Q(sI-A)^{-1}B$.
\end{remark}

To give a concrete construction of the adaptive optimal control, we first introduce the WLS estimator and the random regularization method.

\subsection{The WLS estimation with random regularization} 
To describe the adaptive problem in the standard form, let us denote
\begin{equation}\label{key9}
\theta ^{\mathsf{T}}=[A,B] \text{ and } \varphi(t)=[x^{\mathsf{T}}(t),u^{\mathsf{T}}(t)]^{\mathsf{T}},
\end{equation}
and rewrite the system (\ref{s1}) into:
\begin{equation}
x(t+1)=\theta^{\mathsf{T}}\varphi(t)+Dw(t).
\end{equation}

The  WLS estimates $\{\theta(t),t\ge0\}$ are given by \cite{a1}
	\begin{equation}\label{wls1}
		\theta(t+1)=\theta(t)+a(t)P(t)\varphi(t)[x^{\mathsf{T}}(t)-\varphi^{\mathsf{T}}(t)\theta(t)],
	\end{equation}
	\begin{equation}\label{wls2}
		P(t+1)=P(t)-a(t)P(t)\varphi(t)\varphi^{\mathsf{T}}(t)P(t),
	\end{equation}
	where $P(0)=I_{n+m}$, $\theta^{\mathsf{T}}(0)=[A(0),B(0)]$ are arbitrary deterministic values such that $(A(0),B(0))$ is stabilizable, and
\begin{equation}\label{a1}
	a(t)=1/(\log^2r(t)+\varphi^{\mathsf{T}}(t)P(t)\varphi(t)),
\end{equation}
with
\begin{equation}\label{r}
	r(t)=\|P^{-1}(0)\|+\sum_{k=0}^{t}|\varphi(k)|^2.
\end{equation}
The properties of WLS estimates are presented in the following lemma.
\begin{lemma} \cite{a1}\label{le1}
	Let $\{\theta(t),t \in\mathbb{N}\}$ be defined by (\ref{wls1})-(\ref{r}). Then the
	following properties are satisfied:
\begin{flalign}
		&1) \sup\limits_{t\ge0}\|P^{-\frac{1}{2}}(t)(\theta(t)-\theta)\|<\infty\quad a.s.\notag&&\\
		&2)\sum_{t=0}^{\infty}a(t)\|(\theta(t)-\theta)^{\mathsf{T}}\varphi(t)\|^2<\infty\quad a.s.\notag&&\\
		&3)\lim\limits_{t\to\infty}\theta(t)=\bar{\theta}\quad a.s.\notag&&
\end{flalign}
	where $\bar{\theta}$ is a certain random matrix.
\end{lemma}

From Lemma \ref{le1}, one can see that the WLS estimates are self-convergent, but it should be noted that the family of estimated models may not be stabilizable. In order to obtain stabilizable estimated pair, we adopt the random regularization method introduced in \cite{a1}.

Let us denote the following two matrix functions:
$$Z(s,A,B)\!=\!\left[\left(\begin{array}{cc}
		s_1I_n-A & s_2I_n  \\ 
		-s_2I_n & s_1I_n-A \\ 
	\end{array}\right),\left(\begin{array}{cc}
	B  & 0 \\ 
	0 & B  \\ 
	\end{array}\right)\right],$$
 where $s=s_1+s_2j$ and $j$ is the imaginary unit, and
 \begin{align}
     Y(A,B)&=\prod_{s\in\sigma_{\ge1}(A)}\det\left(Z(s,A,B)Z^{\mathsf{T}}(s,A,B)\right) \notag\\
     &\quad\cdot\prod_{s\in\sigma_{<1}(A)}\det\left(Z(-s,A,B)Z^{\mathsf{T}}(-s,A,B)\right),\notag
 \end{align}
 where $\sigma_{\ge1}\triangleq\{\lambda\in\sigma(A)~\!|~\!|\lambda|\ge1\}$, and $\sigma_{<1}\triangleq\{\lambda\in\sigma(A)~\!|~\!|\lambda|<1\}$.
\begin{lemma}\cite{g2}\label{le2}
    The following statements are equivalent:\\
1) $(A,B)$ is stabilizable;\\
2) $(sI-A,B)$ is of full of row rank for any $s\in\sigma_{\ge1}$;\\
3) $Z(s,A,B)$ is of full of row rank for any $s\in\sigma_{\ge1}$;\\
 4) $Y(A,B)>0$.
\end{lemma}
\vskip 2mm
\begin{definition}[Uniformly stabilizable]\cite{g2}
A family of matrix pairs $\{A(t),B(t); t \in\mathbb{N}\}$ with $A(t)\in \mathbb{R}^{n\times n},B(t) \in \mathbb{R}^{n\times m}$ is said to be uniformly stabilizable if there is a constant $c > 0$ such that
	$$Y(A(t),B(t))\ge c,\text{ for all } t \in\mathbb{N}.$$
\end{definition}
\vskip 2mm

From the boundedness of the matrix sequence $\{P^{-\frac{1}{2}}(t)(\theta-\theta(t)), t\in\mathbb{N}\}$ as established in Lemma \ref{le1}, there exists a bounded sequence $\{\beta^*(t),t\in\mathbb{N}\}$ such that
\begin{equation}\label{b}
    \theta=\theta(t)-P^{\frac{1}{2}}(t)\beta^*(t).
\end{equation}
This serves as a basis for a modification to the WLS estimates:
\begin{equation}\notag
	\theta(t,\beta(t))=\theta(t)-P^{\frac{1}{2}}(t)\beta(t),
\end{equation}
where $\{\beta(t)\in\mathbb{R}^{(m+n)\times n}, t\in\mathbb{N}\} $ is a sequence of bounded matrices to be chosen shortly. For simplicity, we rewrite $\theta(t,\beta(t))$ in the following form:
\begin{equation}\notag
	\theta^{\mathsf{T}}(t,\beta(t))=[\bar{A}(t),\bar{B}(t)].
\end{equation}

Then, we aim to choose suitable $\{\beta(t),t\in\mathbb{N}\}$ to ensure that $\{\bar{A}(t),\bar{B}(t);t\in\mathbb{N}\}$ is uniformly stabilizable.
 From Lemma \ref{le2}, it suffices to select  $\{\beta(t),t\in\mathbb{N}\}$ to guarantee the uniform positivity of $$f_t(\beta(t))\triangleq Y(\bar{A}(t),\bar{B}(t)),$$
 i.e., $f_t(\beta(t)\ge c$ for all $t\in\mathbb{N}$ with some constant $c>0$.
	
	 Next, we give a concrete construction of the sequence $\{\beta(t),t\in\mathbb{N}\}$. Let $\{\eta_t\in\mathbb{R}^{(n+m)\times n}, t\in\mathbb{N}\}$ be a sequence of random variables that are uniformly distributed in the unit ball of $\mathbb{R}^{(m+n)\times n}$ and are also independent of $\{w(t),t\in\mathbb{N}\}$. The sequence $\beta(t)$ is recursively given by the following:
	\begin{align}
		&\beta(0)=0,\notag \\
		&\beta(t)=\begin{cases}
			\eta_t,&\text{if } f_t(\eta_t))\ge\gamma f_t(\beta(t\!-\!1)),\\
			\beta(t \!-\!1),&\text{otherwise},
		\end{cases}\notag
	\end{align}
	where $\gamma\in(1,\sqrt{2})$ is a fixed constant. Thus, a sequence of
	regularized estimates $\{\hat{\theta}(t), t\in\mathbb{N}\}$ can be defined by 
	\begin{equation}\label{wls4}
		\hat{\theta}(t)=\theta(t)-P^{\frac{1}{2}}(t)\beta(t).
	\end{equation}
For the above regularized estimates, we have the following lemma, which can be derived  along the proof line  in \cite{a1}.
	\begin{lemma}\label{le3}
		If Assumption \ref{ass1} is satisfied, then for any admissible control $u(t)$, the family of regularized WLS estimates $\{\hat{\theta}(t),t\in\mathbb{N}\}$ defined by (\ref{wls1})-(\ref{wls4}) has the following properties.
		\par 
		1) Self-convergence, that is, $\hat{\theta}(t)$ converges a.s. to a finite random matrix as $t\to\infty$.
		\par 
		2) The family of matrices $\{A(t),B(t);t\in\mathbb{N}\}$ is uniformly stabilizable where $[A(t),B(t)]=\hat{\theta}^{\mathsf{T}}(t).$
		
		3) 
		$\sum_{s=0}^{t}|(\hat{\theta}(s)-\theta)^{\mathsf{T}}\varphi(s)|^2=o(r(t))+O(1)$, 
		where $r(t)$ is defined by (\ref{r}).

	\end{lemma}
\section{Main Results}\label{section3}
First, we recall the notation
$$\hat{\theta}^{\mathsf{T}}(t)=[A(t),B(t)],$$
and consider the corresponding  ARE:
\begin{equation}\label{are2}
	X(t)=A^{\mathsf{T}}(t)X(t)(I_n+S(t)X(t))^{-1}A(t)+Q,
\end{equation}
where $S(t)=B(t)R^{-1}B^{\mathsf{T}}(t)$. However there is no enough conditions to guarantee that the ARE (\ref{are2}) has a stabilizing solution $X_1(t)$. Fortunately, by Lemma \ref{le3}, $(A(t),B(t))$ is stabilizable, then the modified ARE
\begin{equation}\notag
	X(t)=A^{\mathsf{T}}(t)X(t)(I_n+S(t)X(t))^{-1}A(t)+I,
\end{equation}
has a unique stabilizing solution $X_2(t)$, which can help us to design the adaptive control.
\begin{remark}
    It is worth mentioning that the ARE (\ref{are2}) has at most one stabilizing solution and the stabilizing solution is locally analytic with respect to the matrix parameters (see \cite{b2}). Hence, when $t$ is sufficiently large, the stabilizing solution to the ARE (\ref{are2}) either exists all the time or does not exist at all.
\end{remark}

Following the switching method in \cite{g1}, by certainty equivalence principle, the adaptive control could have the form:
\begin{equation}\label{key4}
	u(t)=L(t)x(t), 
\end{equation}
where $L(t)=-(R+B^{\mathsf{T}}(t)X^*(t)B(t))^{-1}B^{\mathsf{T}}(t)X(t)A(t)$ with
\begin{equation*}
    X^*(t)=\begin{cases}
        X_1(t), \text{ if (\ref{are2}) has a stabilizing solution;}\\
        X_2(t), \text{ otherwise}.
        \end{cases}
\end{equation*}

As shown in \cite{b1}, some diminishing excitation or exploration signals will be helpful for the strong consistency of the estimates but not essentially influence the performance of the system. This inspires us to design the following adaptive control
\begin{equation}\label{c1}
    u^*(t)=L(t)x(t)+v(t)/t^\epsilon, 
\end{equation}
where the sequence $\{v(t),t \in\mathbb{N}\}$  is a white noise independent of $\{w(t),t\in\mathbb{N}\}$ and $\{\eta_t,t\in\mathbb{N}\}$, and $\epsilon\in(0,1/5)$.

\begin{theorem}\label{th1}{\rm (Global Stability)}
    If Assumption \ref{ass1} is satisfied, then under the
adaptive control (\ref{c1}), the system (\ref{s1}) is globally stable in the sense that for any initial state $x(0)$,
\begin{equation}
\limsup\limits_{T\to \infty}\frac{1}{T}\sum_{t=0}^{T}|x(t)|^2<\infty \quad a.s.
\end{equation}
\end{theorem}

Next, we discuss whether the estimates $\hat{\theta}(t)$ given by (\ref{wls4}) converge to the true value $\theta$. 
Let $B_0=[B,D]$ and $\text{Im}(A|B_0)= \text{Im}(B_0)+\cdots+A^{n-1}\text{Im}(B_0)$. Choose $n_1$ linear independent bases $e_1,e_2,\cdots,e_{n_1}\in\text{Im}(A|B_0)$, $n-n_1$ linear independent bases $e_{n_1+1},\cdots,e_{n}\in \text{Im}(A|B_0)^\bot$, and transformation matrix $U=[e_1,e_2,\cdots,e_n]$, then the system (\ref{s1}) can be rewritten into
\begin{equation}\label{s2}
    \bar{x}(t+1)=\bar{A}\bar{x}(t)+\bar{B}u(t)+\bar{D}w(t),
\end{equation}
  where 
  \begin{equation}\label{a}
  \bar{x}(t)=U^{-1}x(t)
  \begin{pmatrix}
      x_1(t)  \\ 
	x_2(t) \\ 
  \end{pmatrix}
,\bar{A}=U^{-1}AU=\begin{pmatrix}
    A_1 & A_2  \\ 
		0 & A_3\\ 
\end{pmatrix},
  \end{equation}
  and
  $$\bar{B}=U^{-1}B=\left(\begin{array}{cc}
		B_1  \\ 
		0 \\ 
	\end{array}\right),\bar{D}=U^{-1}D=\left(\begin{array}{cc}
		D_1  \\ 
		0 \\ 
	\end{array}\right).$$
According to the special form of $\bar{A},\bar{B},\bar{D}$,  the dynamics of  $x_2$ is given by 
 $$x_2(t+1)=A_3x_2(t).$$
 It is obvious that there exists some initial state $x_2(0)$ (e.g., $x_2(0)=0$) such that we cannot obtain the true value $A_3$ no matter what the input control is. In fact, we could obtain the value $(A_1,B_1)$
   for any initial state, which will be shown in the following theorem. But, since $A,B,D$ are unknown, one direct challenge is how to obtain a desired matrix $U$. To address this issue, let us denote a linear space
$$\mathcal{U}(t)=\text{span}\Big\{\eta\Big|\eta\text{ is an eigenvector of }\sum_{s=0}^{t}\varphi(s)\varphi^{\mathsf{T}}(s)$$$$ \quad\quad\quad\text{whose corresponding eigenvalue } \lambda<\log t\Big\}.$$
Then we have the following lemma and its proof can along the proof line of Lemma 2 in \cite{g2}.
 \begin{lemma}\label{le8}
    For the system (\ref{s1}) with the adaptive control (\ref{c1}), there exists some $t_0\ge0$ such that when $t\ge t_0$,  $[\eta_1^{\mathsf{T}},\eta_2^{\mathsf{T}}]^{\mathsf{T}}\in \mathcal{U}(t)$ if and only if  $\eta_1\in \text{Im}(A|B_0)^\bot$ and $\eta_2=0$.
\end{lemma}
\vskip 2mm

From the above lemma, when $t$ is large enough, the space $\text{Im}(A|B_0)^\bot$ can be obtained easily from $\mathcal{U}(t)$. Then, for any transformation matrix $U=[e_1,e_2,\cdots,e_n]$, where $e_1,e_2,\cdots,e_{n_1}$ is $n_1$ bases of $\text{Im}(A|B_0)$ and $e_{n_1+1},\cdots,e_{n}$ is $n-n_1$ bases of $\text{Im}(A|B_0)^\bot$, we have the following theorem.

 \begin{theorem}\label{th2}{\rm (Consistency)}
    Under the same conditions in Theorem \ref{th1}, the parameter estimates $\hat{\theta}^{\mathsf{T}}(t)=[ A(t), B(t)]$ satisfy
     \begin{equation*}
\lim_{t\to\infty } A(t)=U\begin{pmatrix}
A_1 & \bar A_2  \\ 
		0 & \bar A_3
	\end{pmatrix}U^{-1},~~\lim_{t\to\infty } B(t)=B~~~~ a.s.,
     \end{equation*}
      where $\bar A_2,\bar A_3$ are some random matrices which may not equal to the true values $A_2,A_3$ given by (\ref{a}). 
 \end{theorem}
 \vskip2mm
\begin{remark}
    If $(A,B)$ is controllable, then Theorem \ref{th2} implies the strong consistency of $\hat{\theta}(t)$, which have been established in \cite{a1}. However, if $(A,B)$ is only stabilizable, then the estimates of $(A_2,A_3)$ are not guaranteed to be strongly consistent. Fortunately, the above property for estimated models is enough for the optimality  of our adaptive  control (\ref{c1}), which will be shown in the following theorem.
\end{remark}

\begin{theorem}\label{th3}{\rm (Optimality)}
    Under the same conditions in Theorem \ref{th1}, the adaptive control (\ref{c1}) is optimal, that is, 
\begin{equation}
	J(u^*)=\text{tr}(D^{\mathsf{T}}X^*D)\quad a.s.
\end{equation}
where $X^*$ is the stabilizing solution to the ARE (\ref{are1}).
\end{theorem}
\section{Proofs}\label{section4}
\subsection{Proof of Theorem \ref{th1}:}
By Lemma \ref{le3}, there are random matrices $A(\infty)$ and $B(\infty)$ such that
$$\lim\limits_{t\to \infty}A(t)=A(\infty), \quad\lim\limits_{t\to \infty}B(t)=B(\infty)\quad a.s.$$
and that $(A(\infty),B(\infty))$ is stabilizable. For simplicity of the remaining descriptions, we denote
\begin{equation}\label{l2}
\Phi(t)=A(t)+B(t)L(t),
\end{equation}
where $L(t)$  is defined in (\ref{c1}).

 Under adaptive control (\ref{c1}), the system (\ref{s1}) will be
\begin{align}\label{s}
x(t+1)&=Ax(t)+Bu^*(t)+Dw(t)\\
&=\Phi(t)x(t)+\delta(t)+v_1(t)+Dw(t),\notag
\end{align}
where $\delta(t)=(\theta-\hat{\theta}(t))^{\mathsf{T}}\varphi(t)$ and $v_1(t)=Bv(t)/t^\epsilon$.

Since $\Phi(k)$ is stable and convergent, by Cauchy-Schwarz inequality, it is easy to get
\begin{equation}\label{e1}
    |x(t+1)|^2\le m|x(t)|^2+m_1(\delta(t)|^2+|v_1(t)|^2+|Dw(t)|^2),
\end{equation}
where $0<m<1$ and $m_1>0$ are some fixed constants associated with the supremum of the family $\{\Phi(t),t\ge0\}$. Then, it is easy to see that from (\ref{e1}) and Lemma \ref{le3},
\begin{align}
\frac{1}{N}\sum\limits_{t=1}^{N}|x(t)|^2
&=\frac{m}{N}\sum\limits_{t=0}^{N-1}|x(t)|^2+O(\frac{1}{N}\sum\limits_{t=0}^{N-1}|Dw(t)|^2)\notag\\
&\quad +O(\frac{1}{N}\sum_{t=0}^{N-1}|\delta(t)|^2)+O(\frac{1}{N}\sum_{t=0}^{N-1}|v_1(t)|^2)\notag\\
&=\frac{m}{N}\sum\limits_{t=0}^{N-1}|x(t)|^2+O(1)+o(\frac{1}{N}r(N-1)).\notag
\end{align}
Since
\begin{equation}\notag
    r(N)=\|P^{-1}(0)\|+\sum_{t=0}^{N}|\varphi(t)|^2=O(1)+O(\sum\limits_{t=1}^{N}|x(t)|^2),
\end{equation}
it follows that
$$\frac{1}{N}\sum\limits_{t=0}^{N}|x(t)|^2=O(1),$$
which completes the proof.

\subsection{Proof of Theorem \ref{th2}:}
Firstly, we introduce the following lemma whose proof is similar to Lemma 5 in \cite{g2}.
\begin{lemma}\label{le7}
     For the WLS estimator (\ref{wls1})-(\ref{wls4}), let us denote 
     $$\theta=\left(\begin{array}{cc}
		\theta_1  \\ 
		\theta_2 \\ 
	\end{array}\right),\hat{\theta}(t)=\left(\begin{array}{cc}
		\theta_1(t)  \\ 
		\theta_2(t) \\ 
	\end{array}\right),\varphi(t)=\left(\begin{array}{cc}
		\varphi_1(t)  \\ 
		\varphi_2(t) \\ 
	\end{array}\right)$$
where $\theta_1,\theta_1(t),\varphi_1(t)$ have compatible dimensions. If for sufficiently large $t\in\mathbb{N}$ we have 
 $$\frac{\log ^2t\lambda_{\max}(\sum_{s=0}^{t}\varphi_2(s)\varphi_2^{\mathsf{T}}(s))}{ \lambda_{\min}(\sum_{s=0}^{t}\varphi_1(s)\varphi_1^{\mathsf{T}}(s))}=o(1)\quad a.s.$$
 then the following property holds:
$$\lim\limits_{t\to\infty}\theta_1(t)=\theta_1\quad a.s.$$
 \end{lemma}
 
 In order to apply Lemma \ref{le7} to the estimates $\hat{\theta}(t)$ defined in (\ref{wls4}), we proceed to prove that for the system (\ref{s2}) with the adaptive control (\ref{c1}), we have for all sufficiently large $t\in\mathbb{N}$,
 \begin{equation}\label{eq1}
\lambda_{\min}\left(\sum_{s=0}^{t}U_1\varphi(s)\varphi^{\mathsf{T}}(s)U_1^{\mathsf{T}}\right)>t^{1/2}\quad a.s.
 \end{equation}
 \begin{equation}
\lambda_{\max}\left(\sum_{s=0}^{t}U_2\varphi(s)\varphi^{\mathsf{T}}(s)U_2^{\mathsf{T}}\right)=O(1)\quad a.s.
 \end{equation}
where 
$$U_1=\left(\begin{array}{ccc}
	I_{n_1} & 0 &0  \\ 
	0& 0  &I_m\\ 
\end{array}\right)\mathsf{diag}(U^{-1},I_m),$$$$U_2=[0,I_{n-n_1},0]\mathsf{diag}(U^{-1},I_m).$$

For this purpose,  we need to note that the augmented state 
     $$\varphi(t)=[x^{\mathsf{T}}(t),u^{*\mathsf{T}}(t)]^{\mathsf{T}}$$ 
     satisfies the following equation for any $t\in\mathbb{N}$:
\begin{equation}\label{a15}
 \varphi(t+1)=C_t\varphi(t)+G_t\xi(t),
\end{equation}
where
$$C_t=\left(\begin{array}{cc}
	A & B \\ 
	L(t)A & L(t)B  \\ 
\end{array}\right), G_t=\left(\begin{array}{cc}
	D & 0  \\ 
	L(t)D & I_{m}  \\ 
\end{array}\right),$$
$$\xi(t)=(w(t),v(t)/t^\epsilon).$$
Then, we have
\begin{align}\label{cg}
    &\mathsf{diag}(U^{-1},I_m)(C_t,G_t)\mathsf{diag}(U,I_{2m+p})\notag\\
    &=\left(\left(\begin{array}{ccc}
	A_1 & A_2 & B_1   \\ 
	0 & A_3 & 0   \\ 
        * & * & * \\
\end{array}\right),\left(\begin{array}{cc}
	D_1 & 0  \\
        0   & 0  \\
	* & I_m   \\ 
\end{array}\right)\right).
\end{align}

Let us denote
$$H_{t}=\varphi(t)\varphi^{\mathsf{T}}(t),$$
$$O_t=[G_t,C_tG_t,\cdots,C^{n+m-1}_tG_t].$$
It is not hard to verify that
\begin{align}
    &\liminf_{t\to\infty}\frac{1}{t^{3/5}}\lambda_{min}\big(\sum_{s=0}^{t}U_1H_sU^{\mathsf{T}}_1\big)\notag\\
    &\ge\! \liminf_{t\to\infty}\frac{1}{t^{3/5}}\lambda_{min}\big(\sum_{s=0}^{t-1}U_1(C_sH_sC^{\mathsf{T}}_s\!+\!G_s \xi(s)\xi^{\mathsf{T}}(s)G_s^{\mathsf{T}})U^{\mathsf{T}}_1\big)\notag\\
    &\ge \liminf_{t\to\infty}\frac{1}{t^{3/5}}\lambda_{min}\big(\sum_{s=0}^{t-m-n}U_1O_s \xi(s)\xi^{\mathsf{T}}(s)O_s^{\mathsf{T}}U^{\mathsf{T}}_1\big)\notag
\end{align}

Then for any $\eta\in\text{Im}(A|B_0)\oplus\mathbb{R}^n$ with $\|\eta\|=1$, we aim to show that
\begin{equation}\label{n1}
    \liminf_{t\to\infty}\frac{1}{t^{3/5}}\sum_{s=0}^{t}\eta ^{\mathsf{T}}O_s\xi(s) \xi^{\mathsf{T}}(s)O_s^{\mathsf{T}}\eta>0.
\end{equation}
If not, there exists some $\eta^*\in\text{Im}(A|B_0)\oplus\mathbb{R}^n$ with $\|\eta^*\|=1$ such that 
\begin{equation}\label{eq2}
    \liminf_{t\to\infty}\frac{1}{t^{3/5}}\sum_{s=0}^{t}\eta ^{*\mathsf{T}}O_s\xi(s) \xi^{\mathsf{T}}(s)O_s^{\mathsf{T}}\eta^*=0.
\end{equation}
In fact, the equation (\ref{eq2}) means that
\begin{equation}\label{k9}
    \eta^{*\mathsf{T}}[G,CG,\cdots,C^{n+m-1}G]=0,
\end{equation}
where $(C,G)=\lim_{{t}\to\infty}(C_t,G_t)$. Note that $(C,G)$ still has the special form (\ref{cg}), then $\eta^*\in\text{Im}(C|G)$, which is contradictory to (\ref{k9}). Hence, (\ref{n1}) is true, and then (\ref{eq1}) is true.

Next, it is easy to verify that
    \begin{align}
    &\lambda_{\max}\big(U_2(\sum_{s=0}^{t}\varphi(s)\varphi^{\mathsf{T}}(s)  )U_2^{\mathsf{T}}\big)\notag\\
    &=\lambda_{\max}\big(\sum_{s=0}^{t}\mathsf{diag}(0_{n_1},x_2^{\mathsf{T}}(s)x_2(s),0_m) \big)=O(1).\notag
    \end{align}

 Applying Lemma \ref{le7} to $\hat{\theta}^{\mathsf{T}}(t) \mathsf{diag}(U,I_m)$,  
 it follows that
     \begin{align}\label{aa}
         &\lim\limits_{t\to\infty}\hat{\theta}^{\mathsf{T}}(t) \mathsf{diag}(U,I_m)\mathsf{diag}(I_{n_1},0_{n-n_1},I_m)\notag\\
         &=\theta^{\mathsf{T}} \mathsf{diag}(U,I_n)\mathsf{diag}(I_{n_1},0_{n-n_1},I_m)\quad a.s.\notag\\
         &=U\left(\begin{array}{ccc}
		A_1 & 0 & B_1  \\ 
		0 & 0 & 0  \\ 
	\end{array}\right).
     \end{align}
     From (\ref{aa}) and Lemma \ref{le2}, there exist some random matrices $\bar A_2,\bar A_3$ such that
     \begin{equation*}
\lim_{t\to\infty } A(t)=U\begin{pmatrix}
A_1 & \bar A_2  \\ 
		0 & \bar A_3
	\end{pmatrix}U^{-1},~~\lim_{t\to\infty } B(t)=B~~~~ a.s.
     \end{equation*}
     which completes the proof.

\subsection{Proof of Theorem \ref{th3}:}
    By Theorem \ref{th1} and Lemma \ref{le3}, it is obvious that the adaptive control (\ref{c1}) is admissible. Next, we proceed to verify that
\begin{equation}
	J(u^*)=\text{tr}(D^{\mathsf{T}}X^*D).
\end{equation}
Under the adaptive control (\ref{c1}), the system (\ref{s2}) will be
\begin{align}\label{s3}
	x_1(t+1)&=(A_1x_1(t)+A_2x_2(t)+B_1u^*(t))+D_1w(t)\notag\\
	&=(\Phi_1(t)x_1(t)+\Phi_2(t)x_2(t)+v_1(t)+D_1w(t),
\end{align}
and
$$x_2(t+1)=A_3x_2(t),$$
 where $\Phi_1(t)=A_1-B_1L_1(t)$, $\Phi_2(t)=A_2-B_1L_2(t)$ with $[L_1(t),L_2(t)]=L(t)$ and $v_1(t)=B_1v(t)/t^\epsilon$.

Note that the corresponding ARE (\ref{are1}) is equivalent to the following three equations:
\begin{equation}\label{are3}
    X_1=Q_1+A_1^{\mathsf{T}}X_1A_1-A_1^{\mathsf{T}}XB_1(R+B_1^{\mathsf{T}}X_1B_1)^{-1}B_1^{\mathsf{T}}X_1A_1
\end{equation}
$$X_2=Q_2+A_1^{\mathsf{T}}X_1A_2+A_1^{\mathsf{T}}X_2A_3$$
$$-A_1^{\mathsf{T}}XB_1(R+B_1^{\mathsf{T}}X_1B_1)^{-1}(B_1^{\mathsf{T}}X_1A_2+B_1^{\mathsf{T}}X_2A_3)$$
$$X_3=Q_3+A_2^{\mathsf{T}}X_1A_2+A_2^{\mathsf{T}}X_2A_3+A_3^{\mathsf{T}}X_1A_2+A_3^{\mathsf{T}}X_3A_3$$
$$\!-\!(A_2^{\mathsf{T}}X_1B\!+\!A_3^{\mathsf{T}}X_2B_1)(R\!+\!B_1^{\mathsf{T}}X_1B_1)^{-1}(B_1^{\mathsf{T}}X_1A_2\!+\!B_1^{\mathsf{T}}X_2A_3)$$
where
$$\bar{X}\!=\!U^{\mathsf{T}}XU\!=\!\left(\begin{array}{cc}
		X_1 & X_2  \\ 
		X_2^{\mathsf{T}} & X_3 \\ 
	\end{array}\right), \bar{Q}\!=\!U^{\mathsf{T}}XU\!=\!\left(\begin{array}{cc}
		Q_1 & Q_2  \\ 
		Q_2^{\mathsf{T}} & Q_3 \\ 
	\end{array}\right).$$

Since the stabilizing solution to the ARE (\ref{are2}) is locally analytic with respect to the matrix parameters (see \cite{b2}) and some components of the estimates $(A(t),B(t))$ introduced in Theorem \ref{th2} converge to the truth values, without loss of generality, we can assume that ARE (\ref{are2}) admits a stabilizing solution for every $t\in\mathbb{N}$.

From (\ref{s3}), it is not hard to prove that
$$\lim_{t\to\infty}|x_1(t)|^2/t\to0.$$

Let $X_1^*$ be the stabilizing solution to the ARE (\ref{are3}), then we have
\begin{align}
    &\langle X_1^*x_1(t+1),x_1(t+1)\rangle-\langle X_1^*x_1(t),x_1(t)\rangle\notag\\
    &=\langle (\Phi_1^{\mathsf{T}}(t)X_1^*\Phi_1(t)-X_1^*)x_1(t),x_1(t)\rangle\notag\\
    &\quad+\langle X_1^*\Phi_1(t)x_1(t),\Phi_2(t)x_2(t)+v_1(t)+D_1w(t)\rangle\notag\\
    &\quad+\langle X_1^*(\Phi_2(t)x_2(t)+v_1(t)+D_1w(t)),\Phi_2(t)x_2(t)\notag\\
    &\quad+v_1(t)+D_1w(t)\rangle\notag.
\end{align}
Summing above equation from $0$ to $T-1$, we get
\begin{align}\label{e2}
    &\langle X_1^*x_1(T),x_1(T)\rangle-\langle X_1^*x_1(0),x_1(0)\rangle\notag\\
    &=\sum\limits_{t=0}^{T-1}\langle (\Phi_1^{\mathsf{T}}(t)X_1^*\Phi_1(t)-X_1^*)x_1(t),x_1(t)\rangle\notag\\
    &\quad+\sum\limits_{t=0}^{T-1}\langle X_1^*\Phi_1(t)x_1(t),\Phi_2(t)x_2(t)+v_1(t)+D_1w(t)\rangle\notag\\
    &\quad+\sum\limits_{t=0}^{T-1}\langle X_1^*(\Phi_2(t)x_2(t)+v_1(t)+D_1w(t)),\Phi_2(t)x_2(t)\notag\\
    &\quad+v_1(t)+D_1w(t)\rangle.
\end{align}
It is easy to verify that
$$\limsup\limits_{T\to \infty}\frac{1}{T}\sum_{t=0}^{T}x_2^{\mathsf{T}}(t)x_1(t)=\!\limsup\limits_{T\to \infty}\frac{1}{T}\sum_{t=0}^{T} x_2^{\mathsf{T}}(t)v_1(t)\!=\!0,$$
$$\limsup\limits_{T\to \infty}\frac{1}{T}\sum_{t=0}^{T}x_2^{\mathsf{T}}(t)w_1(t)\!=\!\limsup\limits_{T\to \infty}\frac{1}{T}\sum_{t=0}^{T}v_1^{\mathsf{T}}(t)v_1(t)\!=\!0.$$
Combining these with (\ref{e2}), it follows that
\begin{align}
    &\limsup\limits_{T\to \infty}\frac{1}{T}\sum_{t=0}^{T}\langle (\Phi_1^{\mathsf{T}}(t)X_1^*\Phi_1(t)-X_1^*)x_1(t),x_1(t)\rangle\notag\\
    &=\text{tr}(D_1^{\mathsf{T}}X_1^*D_1)=\text{tr}(D^{\mathsf{T}}X^*D).\notag
\end{align}

Then, we have
\begin{align}
J(u^*)
&=\limsup\limits_{T\to \infty}\frac{1}{T}\sum_{t=1}^{T}(x^{\mathsf{T}}(t)Qx(t)+u^{*\mathsf{T}}(t)Ru^*(t))\notag\\
    &=\limsup\limits_{T\to \infty}\frac{1}{T}\sum\limits_{t=0}^{T-1}\langle (Q_1+L^{\mathsf{T}}_1(t)RL_1(t))x_1,x_1(t)\rangle\notag\\
    &=\limsup\limits_{T\to \infty}\frac{1}{T}\sum_{t=0}^{T}\langle (\Phi_1^{\mathsf{T}}(t)X_1^*\Phi_1(t)-X_1^*)x_1(t),x_1(t)\rangle\notag\\
    &=\text{tr}(D^{\mathsf{T}}X^*D).\notag
\end{align}
Hence, Theorem \ref{th3} is true.

\section{Conclusions}\label{section6}
In this paper, we  consider the discrete-time adaptive LQ control problem and provide a complete solution under the condition that the ARE (\ref{are1}) has a stabilizing solution. It is shown that the proposed adaptive control has the  ability to achieve global stability and optimality, even if not all components of the estimated models converge to the true values. As a result of our analysis, we demonstrate that for the adaptive LQ control problem studied in this paper, if the optimal control exists when the parameters are known, it is still possible to design an adaptive optimal control when the parameters are unknown. It is worth mentioning that, similar results as those established in this paper can be extended to the continuous-time case.

\end{document}